\documentclass{amsart}
\usepackage{amsfonts}
\newtheorem{theorem}{Theorem}[section]

\theoremstyle{definition}

\theoremstyle{remark}

\numberwithin{equation}{section}

\newcommand{\C}{\mathbb{C}}
\newcommand{\Z}{\mathbb{Z}}

\begin{document}

\title{Periodic attractors of random truncator maps}

\author{Ted Theodosopoulos}
\address{IKOS Research, 9 Castle Square, Brighton, East Sussex, BN1 1DZ, UK}
\email{ptaetheo@earthlink.net}
\author{Bob Boyer}
\address{Department of Mathematics 
\\ Drexel University \\ Philadelphia, PA 19066}
\email{boyerrp@drexel.edu}
\urladdr{http://madhava.math.drexel.edu/~rboyer/}


\date{June 20, 2006.}



\begin{abstract}
This paper introduces the \textit{truncator} map as a dynamical system on the space of configurations of an interacting particle system.  We represent the symbolic dynamics generated by this system as a non-commutative algebra and classify its periodic orbits using properties of endomorphisms of the resulting algebraic structure.  A stochastic model is constructed on these endomorphisms, which leads to the classification of the distribution of periodic orbits for random truncator maps.  This framework is applied to investigate the periodic transitions of Bornholdt's spin market model.
\end{abstract}

\maketitle

\section{Model Description}
Let $\Omega=[-1,1]^N$ for some positive dimension $N$ and consider a set $\{S_j\}_{j=1}^M$ of mutually exclusive and exhaustive subsets of $\Omega$.  A typical example will be the generalized quadrants, i.e. 
$$S_j = \left\{x \in \Omega \left| sgn (x_i) = \alpha_i \mbox{, } 1 \leq i \leq N \mbox{ and } \sum_{i=0}^{N-1} \alpha_{i+1} 2^i = 2^N+1-2j \right. \right\},$$ 
where the unique set of $\{\alpha_i\}$ denotes the binary decomposition of the integer $j < M$.

Given a mapping $f: \Omega \longrightarrow \Omega$, we define the \textit{truncator} map as the following discrete dynamical system:
\begin{equation}
x(n+1)_i = x(n)_i sgn \left( f \left( x(n) \right)_i \right). \label{eq:trunc1}
\end{equation}
In this paper we specialize to the case of \textit{shuffling} maps, i.e. $f$ which can be expressed as a set of invertible operators $A_j$ associated with each component $S_j$ of $\Omega$.

Specifically, consider the finite group $G=\{1, 2, \ldots, M\}$ endowed with an operation $\circ$ such that, for every $g \in G$, $g \circ g = 1$.  This group is naturally isomorphic to the cyclic product group $\Z_2 \times \Z_2 \times \ldots \times \Z_2$ of $M$ factors, which can be represented as a modulo multiplication group $M_n$ for some large enough $n$ such that $\phi(n)=M$, where $\phi$ is the Euler totient function.  In this setting, assign an orientation reversing invertible $\ell_\infty$ isometry $A_j$ to each component $S_j$ of $\Omega$, with the property that $A_j \left(S_i \right) = S_{i \circ j}$.  The associated \textit{shuffling} map is given by a mapping $\varphi: G \longrightarrow G$ such that $f \left|_{S_i} \right. = A_{\varphi(i)}$.  Using this notation, the resulting truncator dynamics can be described as 
\begin{equation}
x(n+1) = \sum_{i=1}^M A_{\varphi (i)} \left( x(n) \right) {\bf 1}_{S_i} \left(x(n) \right). \label{eq:trunc2}
\end{equation}

These dynamics arise in a variety of settings \cite{palle1,palle2,kawamura1}.  We were driven to study the truncator dynamics because they represent the frozen phase limit ($\beta \rightarrow \infty$) of a class of interacting particle systems describing economic interactions and opinion formation \cite{theo1,theo2,theo3,bornholdt1}.  In this setting, the points $x$ represent configurations of a spin network and the shuffling map represents the interaction Hamiltonian that describes the influence of local and global effects to the flipping of individual spins.

Another setting where such truncator dynamics arise is that of random Boolean networks.  Often such models are used to describe regulatory networks (e.g. genetic or metabolic networks in biology \cite{troein1,troein2,troein3}) and they are also used to describe instances of the satisfiability problem \cite{parisi}.  In this latter setting, global optimization algorithms are constructed to flip the values of Boolean variables populating the nodes of a graph in such a way as to maximize the probability that the clauses represented by the graph connections are simultaneously satisfied.  

Our goal in this paper is to characterize the periodic attractors of the truncator map.  Specifically we consider random endomorphisms of $G$ \cite{diaconis,peres} and derive the distribution of periodic orbits of the resulting random truncator dynamics.  Of course the full truncator map (\ref{eq:trunc1}) is generically chaotic \cite{gilmore}, because there is sensitivity to initial conditions in the neighborhood of the boundaries between the components $S_j$ (e.g. the axes, when the components are generalized quadrants).  Here we will restrict our attention to shuffling maps and the resulting restricted truncator dynamics (\ref{eq:trunc2}) which captures the spectrum of periodic attractors.  In a later step we plan to use this analysis as a building block for understanding the transitions between the basins of attraction of the periodic attractors we describe here.

\section{Algebraic Dynamics}

In order to better describe the orbits of (\ref{eq:trunc2}) we define a new, noncommutative operation on $G$.  This operation encodes the action of the shuffling map $\varphi$ on $G$:
$$g_1 \ast g_2 = g_1 \circ \varphi (g_2).$$
Abusing notation and identifying each $x \in \Omega$ with the index of the component $S_i$ in which it lies (i.e. the $i$ such that ${\bf 1}_{S_i} (x) = 1$), and subsequently every index with the corresponding member of $G$, we can describe every orbit of (\ref{eq:trunc2}) as a sequence:
\begin{equation}
g \rightarrow g^{\ast 2} \doteq g \ast g \rightarrow g^{\ast 3} \doteq (g \ast g) \ast (g \ast g) \rightarrow \cdots.  \label{eq:path}
\end{equation}

Here is a list of some preliminary results for this algebraic structure:

\begin{theorem}
\label{homo1}
If $\varphi$ is a homomorphism with respect to $\circ$ then it is also a homomorphism with respect to $\ast$.  Conversely, if $\varphi$ is a surjective homomorphism with respect to $\ast$ then it is also a homomorphism with respect to $\circ$.
\end{theorem}
\begin{proof}
The first statement follows immediately from the definition of the $\ast$ operation, since for any $g_1, g_2 \in G$, both $\varphi(g_1 \ast g_2)$ and $\varphi(g_1) \ast \varphi(g_2)$ are equal to $\varphi(g_1) \circ \varphi^{(2)} (g_2)$ (where $\varphi^{(k)}$ denotes the $k$-fold iteration of $\varphi$).  For the second statement, we observe that, for every $g_2 \in {\rm Im} \varphi$, there exists some $g_3 \in G$ such that $g_2 = \varphi (g_3)$ and thus, $$\varphi(g_1 \circ g_2) = \varphi \left(g_1 \circ \varphi(g_3) \right) = \varphi(g_1 \ast g_3) = \varphi(g_1) \ast \varphi(g_3) = \varphi(g_1) \circ \varphi^{(2)} (g_3) = \varphi(g_1) \circ \varphi(g_2).$$
\end{proof}

\begin{theorem}
\label{fixedpoint}
For any $\varphi$ (not necessarily a homomorphism), if the $\ast$ operation is commutative, then $\# \varphi^{-1} (1) = 1$ and $\varphi^{-1} (1) = 1$ is the unique attractor of (\ref{eq:trunc2}).
\end{theorem}
\begin{proof}
The assumption implies that for every $g_1, g_2 \in G$,
\begin{eqnarray*}
g_1 \ast g_2  = g_2 \ast g_1 & \Longleftrightarrow & g_1 \circ \varphi(g_2) = g_2 \circ \varphi(g_1) \\
& \Longleftrightarrow & g_1 \circ \varphi(g_1) = g_2 \circ \varphi(g_2) \\
& \Longleftrightarrow & g_1 \ast g_1  = g_2 \ast g_2.
\end{eqnarray*}
But this implies that there exists a unique $\hat{g} \in G$ such that all points $x \in \Omega$ move into $S_{\hat{g}}$ after one step of (\ref{eq:trunc2}).  Consider $\hat{g}$ itself.  Since it remains fixed, it is a fixed point of (\ref{eq:trunc2}).  This implies that $\varphi(\hat{g})=1$.  Any $h \in \varphi^{-1} (1)$ is a fixed point since $h \ast h = h \circ \varphi(h) = h$.  But if there was any other member of $\varphi^{-1} (1)$ different from $\hat{g}$, it would have to move to $\hat{g}$ in one step as we have already seen, refuting its stationarity.  No other attractors are possible since all points converge to $\hat{g}$ in one step.  Therefore, $\hat{g}$ is the unique fixed point of (\ref{eq:trunc2}).
\end{proof}

As an example, consider the case $N=2$ and the map $f(z) = {\frac {z + z^{-1}}{2}}$, where we think of $\Omega$ as the unit $\ell_\infty$ ball in ${\C}$.  Of course $f^{-1} (z) = z - \sqrt{z^2 -1}$ which leads us to conclude that this is indeed a shuffling map, with $A_1 (z) = z$, $A_2 (z) = \bar{z} e^{i \pi}$, $A_3 (z) = z e^{i \pi}$ and $A_4 (z) = \bar{z}$, and therefore $\varphi(1) = 4$, $\varphi(2) = 3$, $\varphi(3) = 2$ and $\varphi(4) =1$.  We observe that, with this choice of $\varphi$, $1 \ast 1 = 2 \ast 2 = 3 \ast 3 = 4 \ast 4 = 4 = \varphi^{-1} (1)$.  This, according to Theorem \ref{fixedpoint}, every $z \in \Omega$ with positive real and negative imaginary parts will be a fixed point for the following dynamics
\begin{eqnarray*}
{\rm Re} z(n+1) & = & {\rm Re} z(n) sgn \left({\rm Re} \left ({\frac {z + z^{-1}}{2}} \right) \right) \\
{\rm Im} z(n+1) & = & {\rm Im} z(n) sgn \left({\rm Im} \left ({\frac {z + z^{-1}}{2}} \right) \right).
\end{eqnarray*}

Now, for every $\varphi: G \rightarrow G$ consider a new multiplication in $G$ defined so that it satisfies
$$g \otimes g = \varphi (g)$$
and so that it is left-distributive\footnote{One can easily check that, if $1 \in \ker \varphi$, then $\otimes$ is bilaterally distributive.} with respect to the addition defined by $\circ$.  Note that $\otimes$ is not necessarily associative, e.g. $\{(g \otimes g) \otimes g\} \circ \{g \otimes (g \otimes g) \} = [g, g^{\otimes 2}] = \varphi \left(g \circ \varphi(g) \right) \circ \varphi(g) \circ \varphi^{(2)} (g)$ which can be different from $1$ when $\varphi$ is not a homomorphism.  Due to this potential non-associativity, we must be careful about defining $\otimes$ powers.  In particular, let $\alpha_\ell \alpha_{\ell-1} \cdots \alpha_2 \alpha_1$ be the binary decomposition of the integer $k>1$.  Then define $g^{\otimes k} = \left( \alpha_\ell g^{\otimes 2^{\ell-1}} \right) \otimes \left( \alpha_{\ell -1} g^{\otimes 2^{\ell-2}} \right) \otimes \cdots \otimes \left( \alpha_2 g^{\otimes 2} \right) \otimes (\alpha_1 g)$, where $g^{\otimes 2^j}= \left(g^{\otimes 2^{j-1}} \right) \otimes \left(g^{\otimes 2^{j-1}} \right)$.

If $\varphi$ is a homomorphism, then $\otimes$ is commutative because for any $g, h \in G$,
$$(g \circ h) \otimes ( g \circ h) = \varphi (g \circ h ) = \varphi (g) \circ \varphi (h)$$
while, the distributive property implies that
$$(g \circ h) \otimes ( g \circ h) = (g \otimes g) \circ (h \otimes h) \circ (g \otimes h) \circ (h \otimes g) = \varphi (g) \circ \varphi (h) \circ (g \otimes h) \circ (h \otimes g)$$
and therefore $(g \otimes h) \circ (h \otimes g)=1$ which implies commutativity.  For the same reason, a general $\varphi$ leads to the following identity:
$$(g \otimes h) \circ (h \otimes g) = \varphi (g \circ h) \circ \varphi (g) \circ \varphi (h).$$
We define the commutator of two elements $g,h$ of $G$ as $[g,h] \doteq (g \otimes h) \circ (h \otimes g)$ and say that $g$ commutes with $h$ if $[g,h] = 0$.  Also observe that $1$, the identity of the addition $\circ$, is a trapping element of $G$ with respect to $\otimes$ since for any $g \in G$,
$$\varphi(g)= g \otimes g = (g \circ 1) \otimes g = (g \otimes g) \circ (1 \otimes g) = \varphi (g) \circ (1 \otimes g)$$
and therefore $1 \otimes g = g \otimes 1 = 1$.

\section{Polynomial Roots}

Now consider the ring of polynomials in $G$ using this multiplication and coefficients from $\Z_2$.  The action of $\Z_2$ on $G$ is modeled as exterior multiplication of $\Z_2$ on $M_n$, the modulo multiplication group that represents $G$, i.e. $0 \cdot g = 1$ and $1 \cdot g = g$.

\begin{theorem}
\label{polynomial}
For any $\varphi \in \hom (G,\circ)$, $g \in G$ and $p>0$,
\begin{equation}
g^{\ast p} = \bigodot_{k=0}^{p-1} \gamma_{k,p} \varphi^{(k)} (g), \label{eq:poly}
\end{equation}
where $\varphi^{(0)} (g) = g$,
\begin{equation}
\gamma_{k,p} \equiv \left( \gamma_{k,p-1} + \gamma_{k-1,p-1} \right) \pmod{2}, \label{eq:pascal}
\end{equation}
and $\gamma_{0,0}=1$, $\gamma_{k,0}=0$ for $k<0$.
\end{theorem}
\begin{proof}
Observe that the expression for $g^{\ast n}$ is a polynomial in $G$ as described above, of degree $2^k$, since $\varphi^{(k)} (g) = g^{\otimes 2^k}$, as can be easily checked using induction.  Observe further that the coefficients of these polynomials obey the binary version of the Pascal triangle.  We will show this using induction in $p$.  The desired result clearly holds for $p=1$ since $g = \gamma_{0,1} g$ and $\gamma_{0,1}= 1$.  Assume the desired result holds for $p$.  Then
\begin{eqnarray*}
g^{\ast (p+1)} & = & g^{\ast p} \ast g^{\ast p} = g^{\ast p} \circ \varphi (g^{\ast p}) = \left\{ \bigodot_{k=0}^{p-1} \gamma_{k,p} \varphi^{(k)} (g) \right\} \circ \varphi \left( \bigodot_{k=0}^{p-1} \gamma_{k,p} \varphi^{(k)} (g) \right) = \\ 
& = & \left\{ \bigodot_{k=0}^{p-1} \gamma_{k,p} \varphi^{(k)} (g) \right\} \circ \left\{ \bigodot_{k=1}^{p} \gamma_{k-1,p} \varphi^{(k)} (g) \right\} = \\
& = & \left(\gamma_{0,p} g \right) \left(\gamma_{p-1,p} \varphi^{(p)} (g) \right) \bigodot_{k=1}^{p-1} \left((\gamma_{k-1,p} + \gamma_{k,p}) \pmod{2} \right) \varphi^{(k)} (g)
\end{eqnarray*} 
because when both $\gamma_{k,p}$ and $\gamma_{k-1,p}$ are both equal to $1$, then the corresponding term contains $\varphi^{(k)} \circ \varphi^{(k)}$ and therefore vanishes.
\end{proof}

Let's define the period of an element $g$ as
\begin{equation}
p^\ast (g) = \min \{i>1 | g^{\ast i} = g\} -1   \label{eq:periodef}
\end{equation}
where we understand the minimum of any empty set to be equal to $\infty$.  Using this concept we can summarize a set of necessary and sufficient conditions for a truncator map to possess limit cycles of particular periods as follows:

\begin{theorem}
\begin{enumerate}
\item For a general $\varphi$, $g \in \ker \varphi \Longleftrightarrow p^\ast (g) = 1$.
\item If $g$ commutes with $g^{\otimes 2}$ and $1 \in \ker \varphi$, $g \in \ker \varphi^{(2)} \setminus \ker \varphi \Longleftrightarrow p^\ast (g) = 2$.
\item Let $\bigtriangleup \doteq \{g \in G | g^{\ast 2} = g^{\otimes 4} \}$.  If $g$ commutes with $g^{\otimes 2}$ and $g^{\otimes 4}$ and $1 \in \ker \varphi$, 
$$g \in \varphi^{-1} ({\rm Im}\varphi \cap \bigtriangleup ) \setminus \ker \varphi \Longleftrightarrow p^\ast (g) = 3.$$
\end{enumerate}
\end{theorem}
\begin{proof}
Using (\ref{eq:periodef}) we see that $p^\ast (g)$ is equal to $1$ iff $g^{\ast 2} = g$ which is true iff $g \in \ker \varphi$, thus proving the first statement of the theorem.  On the other hand, we clearly have $g^{\ast 3} = g \circ g^{\otimes 2} \left(g \circ g^{\otimes 2} \right)^{\otimes 2} = g \circ g^{\otimes 4} \circ \left[g, g^{\otimes 2} \right]$.  But when $[g, g^{\otimes 2}] = 1$, $p^\ast (g) \leq 2$ iff $g \in \ker \varphi^{(2)}$.  Together with the previous statement, we have proved the second statement of the theorem.
In this case we have 
\begin{equation}
g^{\ast 4} = g \circ g^{\otimes 2} \circ g^{\otimes 4} \circ g^{\otimes 8} \circ [g, g^{\otimes 2}] \circ [g, g^{\otimes 2}]]^{\otimes 2} \circ [g, g^{\otimes 4}] \circ \left[g \circ g^{\otimes 4}, [g, g^{\otimes 2}] \right].  \label{eq:gast4}
\end{equation}
When $g$ commutes with $g^{\otimes 2}$ and  $g^{\otimes 4}$, (\ref{eq:gast4}) simplifies to $g^{\ast 4} = g \circ g^{\otimes 2} \circ g^{\otimes 4} \circ g^{\otimes 8}$.  Since $\varphi(g) \in \bigtriangleup$, $g^{\otimes 2} \circ g^{\otimes 4} \circ g^{\otimes 8} = 1$ and therefore $g^{\ast 4} = g$ which implies $p^\ast (g) \leq 3$.  Requiring that $g \not \in \ker \varphi$ is sufficient to complete the proof of the last statement in the theorem because $g \in \ker \varphi^{(2)} \setminus \ker \varphi \Longrightarrow g^{\otimes 2} \circ g^{\otimes 4} \circ g^{\otimes 8} = g^{\otimes 2} \neq 1$.
\end{proof}

\section{Random Maps}

Let $\mu \in {\mathcal M}_1 \left( G^G \right)$ be a probability measure on the set of maps from $G$ to itself.  This can be described as a sequence of measures $\nu_g \in {\mathcal M}_1 (G)$ on $G$ indexed by the elements of $G$, such that for every $g, h \in G$:
$$\nu_g (h) = \mu \left( \varphi (g) = h \right).$$
A uniformly random $\varphi$ maps each element $g$ to $1$ with probability $M^{-1}$.  Thus the number of $g$ that are mapped to $1$ is a binomial random variable:
$$\lambda \left( \left|\ker \varphi \right| = k \right) = \left( \begin{array}{c} M \\ k \end{array} \right) M^{-M} (M-1)^{M-k},$$
and therefore,
\begin{theorem}
Let $\lambda$ be the uniform measure in ${\mathcal M}_1 \left( G^G \right)$.  Then:
$$\lim_{M \rightarrow \infty} \lambda \left( \left|\ker \varphi \right| = k \right) = (ek!)^{-1}.$$
\end{theorem}
We proceed by defining a transition matrix $\Phi$ such that for every pair $(i,j) \in G^2$, $\Phi_{i,j} = \mu(i^{\ast 2}=j)$.  We consider a stochastic process on $G$ which propagates according to (\ref{eq:path}) with iid choices of $\varphi$ in every draw.  Observe that
$$\mu \left(i^{\ast 3} =j \right) = \sum_{k \in G} \mu \left(i^{\ast 2} =k \right) \mu \left(k^{\ast 2} =j \right) = \sum_{k \in G} \Phi_{ik} \Phi_{kj} = \left(\Phi^2 \right)_{ij}.$$
Now let $\imath$ be the identity mapping on $G$, and define an addition $+$ in $G^G$ such that $(\varphi_1 +\varphi_2)(g) = \varphi_1(g) \circ \varphi_2(g)$.  Then, we can express $\ast$ powers of $g$ as
$$g^{\ast p} = (\varphi + \imath)^{(p-1)} (g).$$
Then
$$g^{\ast p} = g \Longleftrightarrow \varphi \left(\bigodot_{k=0}^{p-2} (\varphi + \imath)^{(k)} (g) \right) = 1.$$
Thus, when it is finite, $p^\ast$ can be computed as the first passage time into $\ker \varphi$ of a Markov chain on $G$ with transition matrix elements:
\begin{eqnarray*}
{\rm Pr} \left( \bigodot_{k=0}^p \left(\varphi + \imath \right)^{(k)} (g) = j \left| \bigodot_{k=0}^{p-1} \left(\varphi + \imath \right)^{(k)} (g) = i  \right. \right) & = & \mu \left(\left(\varphi + \imath \right)^{(p)} (g) = i \circ j \right) \\
& = & \left(\Phi^p \right)_{g, i \circ j}.
\end{eqnarray*}
On the other hand when $p^\ast (g) = \infty$ ($g$ is transient for the dynamics), the Markov chain never enters $\ker \varphi$.  This is a time inhomogeneous process as seen in the expression for the $p$-step transition probabilities:
\begin{eqnarray*}
& & {\rm Pr} \left( \bigodot_{k=0}^p \left( \varphi + \imath \right)^{(k)} (g) =j \left| g=i \right. \right) = \\
& = & \sum_{k_1,k_2,\ldots,k_{p=2} \in G} \Phi_{i , i \circ k_1} \left(\Phi^2 \right)_{i \circ k_1, i \circ k_1 \circ k_2} \cdots \left(\Phi^p \right)_{\bigodot_{j=1}^{p-1} k_j \circ i, \bigodot_{j=1}^p k_j \circ i }.
\end{eqnarray*}

\section{Synchronous Spin Market Dynamics}

In this section we show how to map a spin model of market microstructure onto the class of truncator dynamics.  The state space $X$ of the model we want to consider is the set of spin configurations on a lattice on the $d$-dimensional torus\footnote{Here we use the notation ${\mathcal T}^d$ to denote the object $\underbrace{ {\mathcal S}^1 \times \ldots \times {\mathcal S}^1}_d$.} $Y \doteq \left({\mathcal Z}/L \right)^d \subset {\mathcal T}^d$, i.e. $X \subset \{-1,1 \}^Y$, for an appropriately chosen $L$ so that $|Y|=N$.  The path of a typical element of $X$ is given by $\eta: Y \times \aleph \longrightarrow \{-1,1\}$ and each site $x \in Y$ is endowed with a (typically $\ell_1$) neighborhood ${\mathcal N} (x) \subset Y$ it inherits from the natural topology on the torus ${\mathcal T}^d$. 

We construct a discrete time Markov process with synchronous transitions updating all the spins simultaneously.  We proceed to construct a transition matrix for the spins, based on the following interaction potential:
$$h(x,n) = \sum_{y \in {\mathcal N}(x)} \eta(y,n) - \alpha \eta(x,n) N^{-1} \left|\sum_{y \in Y} \eta(y,n) \right|,$$
where $\alpha>0$ is the coupling constant between local and global interactions.  At time $n$ the spins change to $+1$ with probability $p^+ \doteq \left( 1+ \exp \left\{- 2\beta h \left(x,n \right) \right\} \right)^{-1}$ and to $-1$ with probability $p^- = 1- p^+$, where $\beta$ is the normalized inverse temperature.

Let $f: X \longrightarrow X$ be such that for all $x \in Y$, $\left(f(\eta( \cdot, n)) \right) (x) \doteq \eta(x,n) h(x,n) = \sum_{y \in {\mathcal N}(x)} \eta(y,n) \eta(x,n) - \alpha N^{-1} \left|\sum_{y \in Y} \eta(y,n) \right|$.  It is easy to check that the frozen phase of this system ($\beta \rightarrow \infty$) is a shuffling truncator map, as described above in (\ref{eq:trunc2}).  In the frozen phase, the transitions are deterministic (each row of $\Phi$ has only one nonzero element).  High but finite values of $\beta$ lead to the introduction of some genuine randomness in the transition matrix.

To illustrate this procedure, let's consider the above spin market model with $N=4$ and $d=1$ and standard nearest neighbor topology in ${\mathcal S}^1$.  Consider first $g=16$ which represents the quadrant $(-1,-1,-1,-1)$.  We realize there are two cases.  When $\alpha <2$, $\varphi(g) = 1$, which represents the quadrant $(1,1,1,1)$; thus when $\alpha<2$ (subcritical regime), $g=16 \in \ker \varphi$ and therefore $p^\ast (16) = 1$.  Notice that by symmetry, the same is true for $g=1$.  On the other hand when $\alpha >2$ (supercritical regime), $\varphi (16) = 16$ and therefore $16 \ast 16 =1$ and $\varphi(1) = 16$.  Thus, $p^\ast (16) = p^\ast (1) = 2$.

Next consider the element $g=15$ representing the quadrant $(-1,-1,-1,1)$.  Once again there are two cases, but this time they are separated by $4$ rather than $2$.  Specifically, when $\alpha<4$, $\varphi (15) = 12$, representing the quadrant $(-1,1,-1,-1)$, and thus $15 \ast 15 = 6$, representing quadrant $(1,-1,1,-1)$.  Proceeding from $g^{\ast 2} = 6$ we see that $\varphi (6) = 16$ and therefore $6 \ast 6 = 11$, representing quadrant $(-1,1,-1,1)$.  Thus, when $\alpha<4$ (the relevant subcritical regime), $p^\ast (15) = \infty$, draining into the period 2 attractor $6 \rightarrow 11 \rightarrow 6 \rightarrow \cdots$.

On the other hand when $\alpha>4$ (the relevant supercritical regime) $\varphi (15) = 16$ and thus $15 \ast 15 = 2$, representing quadrant $(1,1,1,-1)$.  Continuing from $g^{\ast 2} = 2$ we see that $\varphi (2) = 16$ and therefore $2 \ast 2 = 15$.  So we conclude that in the supercritical regime, $p^\ast (15) = 2$.

\section{Conclusions and Next Steps}

We have presented a new methodological framework for analyzing a class of random symbolic dynamics.  This framework draws on the iterated function systems (IFS) literature to identify Boolean maps with Boolean expressions, thus constructing an algebraic structure akin to the modulo multiplication groups.  This structure in turn helps clarify the qualitative properties of the underlying interaction Hamiltonian by exhibiting parameter ranges which lead to different algebraic properties.

We have shown constructively that large classes of symbolic dynamics, including random Boolean networks, can be described in terms of our proposed truncator maps.  In the case of the Bornholdt spin market microstructure model we have shown examples of fixed points, period 2 cycles as well as transient points in configuration space.  We proceeded to show that non-zero temperature can be accommodated by constructing a Markov chain in the space of automorphisms of our ring structure.  In a particularly simple case we were able to compute explicitly the thermodynamic limit of the number of fixed points.  This analytical result lends support to the conjecture that as the number of agents increase, with overwhelming probability there are but very few fixed points.

A natural next step is to extend the analysis presented here beyond shuffling maps to general truncator dynamics.  Such an extension will involve long memory as the iterated images become intertwined.  On the other hand the resulting global mixing is likely to induce ergodic properties missing in the case of pure shuffling.

Furthermore, even in the case of shuffling maps, the solution of the inhomogeneous exit problem identified above as a way to represent the spectrum of the truncator dynamics remains generally open.  We plan to address this problem explicitly in future work.

\bibliographystyle{amsalpha}

\end{document}